\documentclass[a4paper]{amsart}
\usepackage{amssymb}
\usepackage{xy}
\let\<\langle
\let\>\rangle
\DeclareMathOperator{\Id}{Id}
\DeclareMathOperator{\pr}{pr}
\renewcommand{\MR}[1]{}

\numberwithin{equation}{section}

\theoremstyle{plain}
\newtheorem{theorem}[equation]{Theorem}
\newtheorem{proposition}[equation]{Proposition}
\newtheorem{lemma}[equation]{Lemma}
\newtheorem{corollary}[equation]{Corollary}

\theoremstyle{definition}
\newtheorem{definition}[equation]{Definition}
\newtheorem{example}[equation]{Example}

\theoremstyle{remark}
\newtheorem{remark}[equation]{Remark}

\title{Natural Frobenius Submanifolds}

\begin{document}
\author{Jiezhu Lin}
\address{Department of Mathematics, Sun Yat-Sen University, Guangzhou 510275, China}
\email{jiezhu.lin@ens.fr}

\thanks{The research was partially supported by
China-France-Russian mathematics collaboration grant, No.
34000-3275100, from Sun Yat-Sen University.}

\begin{abstract}
I.A.B. Strachan introduced the notion of a natural Frobenius
submanifold of a Frobenius manifold and gave a sufficient but not
necessary condition for a submanifold to be a natural Frobenius
submanifold. This paper will give a necessary and sufficient
condition and classify the natural Frobenius hypersurfaces.
\end{abstract}



\maketitle

\setcounter{section}{-1}

\section{Introduction}\label{section0}

\subsection{Saito structure and Frobenius manifold structure}
\label{subsection0a} Frobenius manifolds were introduced and
investigated by B. Dubrovin as the axiomatization of a part of the
rich mathematical structure of the Topological Field Theory (TFT):
cf. \cite{D, Hert, Mani}

A Frobenius manifold (or called Frobenius structure on $M$) is a
quadruple $(M, \circ, \overline g, e, \mathcal{E})$. Here $M$ is a
manifold in one of the standard categories ($C^\infty$, analytic,
...), $\overline g$ is a metric on $M$ (that is, a symmetric,
non-degenerate bilinear form, also denoted by $\<, \>$), $\circ$
is a commutative and associative product on $TM$ and depends
smoothly on $M$, such that if $\overline \nabla$ denote the
Levi-Civita connection of $\overline g$, then

a) $\overline \nabla$ is flat;

b) $\overline g( X \circ Y, Z) =\overline g( X, Y \circ Z)$, for
any $X, Y, Z \in TM$;

c) the unit vector field e is covariant constant w.r.t. $\overline
\nabla$
\begin{eqnarray*}
\overline \nabla e = 0;
\end{eqnarray*}

d) let
\begin{equation*}
c(X, Y, Z):=\overline g(X \circ Y, Z)
\end{equation*}
(a symmetric 3-tensor). We require the 4-tensor
\begin{equation*}
(\overline \nabla _Z c)(U,V,W)
\end{equation*}
to be symmetric in the four vector fields $U,V,W, Z$.

e) A vector field $\mathcal{E}$ must be determined on $M$ such
that
\begin{align}
\overline \nabla (\overline \nabla \mathcal{E}) &= 0;\label{nablanablaE}\\
\mathcal{L}_{\mathcal{E}}(\circ)&=\circ;\label{LEcirc}\\
\exists D\in\mathbb{C},\quad\mathcal{L}_{\mathcal{E}} (\overline g)&=D \cdot \overline g.\label{LEg}
\end{align}

\begin{remark}\label{remark SABB}
In this definition, because the metric $\overline g$ is flat and
the unit field $e$ is covariant constant w.r.t. $\overline \nabla$,
then \eqref{LEg} implies \eqref{nablanablaE}.

Good reference is the last chapter in \cite{Sabb}.
\end{remark}

There are several equivalent ways to describe a Frobenius
structure. One way, called Saito structure, is recalled here:

\begin{definition}\label{deifinition: Saito structure}
Let $M$ be a complex analytic manifold of dimension~$m$. A Saito
structure on $M$ (without metric) consists of the following data:

1) a flat torsion free connection $\overline \nabla$ on the tangent
bundle $TM$;

2) a symmetric Higgs field $\Phi$ on the tangent bundle $TM$, that
is, $\Phi$ is an $\mathcal{O}_M-$linear map $\Phi$:
$\mathcal{O}(TM) \rightarrow \Omega_M^{1} \otimes \mathcal{O}(TM)$
such that
\begin{equation*}
\Phi_X \Phi_Y = \Phi_Y \Phi_X;
\end{equation*}
3) two global sections (vector fields) $e$ and $\mathcal{E}$ of $\Theta_M$, respectively called unit field and Euler field of the structure.

These data are subject to the following conditions:

a) the meromorphic connection $\widetilde{\overline{\nabla}}$ on
the bundle $\pi^*TM$ on $\mathbb{P}^1 \times M$ defined by the
formula
\begin{eqnarray*}
\widetilde{\overline{\nabla}}=\pi^* {\overline \nabla} + \frac{\pi^* \Phi}{z} - (\frac{\Phi(\mathcal{E})}{z} + \overline \nabla \mathcal{E})\frac{dz}{z}
\end{eqnarray*}
is integrable;

b) the field $e$ is $\overline \nabla$-horizontal (i.e., $\overline
\nabla e = 0$) and satisfies $\Phi_e = - \Id$ (i.e., the product
$\circ$ associated to $\Phi$ has $e$ as a unit field).
\end{definition}

\begin{definition}\label{deifinition: Saito structure with metric}
Let $M$ be a complex analytic manifold of dimension~$m$. A Saito
structure on $M$ with metric consists of a Saito structure
$(\overline \nabla, \Phi, e, \mathcal{E})$ and of a metric
$\overline g$ on the tangent bundle, satisfying the following
properties:

(1) $ \overline \nabla \overline g = 0$ (hence $\overline \nabla$
is the Levi-Civita connection of $\overline g$);

(2) $\Phi^*=\Phi$, i.e., for any local section $X$ of $\Theta_M$, $\Phi^*_X=\Phi_X$, where $^*$ denotes the adjoint w.r.t. $\overline g$;

(3) there exists a complex number $D \in \mathbb{C}$ such that
\begin{equation*}
\overline{\nabla}\mathcal{E}+(\overline{\nabla}\mathcal{E})^*=D\cdot \Id;
\end{equation*}
\end{definition}

\begin{proposition}[\cite{D,Sabb}]\label{proposition3}
On any manifold $M$, there is an equivalence between a Saito structure with metric and a Frobenius structure.
\end{proposition}

\subsection{Frobenius submanifolds}
\label{subsection0b} In \cite{Stra1} the author considers
Frobenius structures defined on open subsets of $\mathbb{R}^n$ or
$\mathbb{C}^n$ and their (natural) Frobenius submanifolds. In
\cite{Stra2}, the author studied the submanifolds $N$ of a semi-simple
Frobenius manifold $M$ with the Euler vector field $\mathcal{E}$
tangent to $N$. We now generalize the definition of a natural
Frobenius submanifold for any Frobenius manifold in the following
way:

Let $(M,\overline{g}, \circ,e, \mathcal{E})$ be a Frobenius
manifold, where $e$ is the unit vector field, $\mathcal{E}$ is the
Euler vector field. Let $N$ be a submanifold of $M$ such that the
metric $\overline{g}$ restricted to $N$, denoted by $g$, is
non-degenerate. So for any tangent vector fields $X,Y\in
\Gamma(U,TN)$ we can define a new product in $TN$ by $X*
Y:=\pr(X\circ Y)$. Similarly we set
$e_{N}:=\pr(e),\mathcal{E}_{N}:=\pr(\mathcal{E})$, where $\pr:
TM\longrightarrow TN$ is the orthogonal projection on $TN$ w.r.t.
$\overline g$. We set $TN^\bot = \lbrace \xi \in TM \mid
\forall X \in TN,\, \<X, \xi \>=0 \rbrace$. So for any vector field
$X \in TM$, we have the decomposition:
\begin{equation*}
X = \pr(X)+X^\bot,
\end{equation*}
where $\pr(X) \in TN$, $X^\bot \in TN^\bot$.

\begin{definition}\label{deifinition: Frobenius submanifold}
The submanifold $N$ is called a \emph{Frobenius submanifold} of
the Frobenius manifold $(M,\overline{g},\circ,e, \mathcal{E})$ if
the induced structure $(N,g,\ast,e_{N}, \mathcal{E}_{N})$ on $N$
is a Frobenius manifold structure.
\end{definition}

\begin{definition}\label{deifinition: Natural Frobenius submanifold}
The Frobenius submanifold $N$ of $(M,\overline{g},\circ,e,
\mathcal{E})$ is called \emph{natural} if $TN$ is left invariant
by the product $\circ$.
\end{definition}

In \cite{Stra1} the author gave a sufficient condition for a
submanifold $N$ to be a natural Frobenius submanifold:

\begin{theorem}[\cite{Stra1}]\label{Theorem4}
Let $N$ be a flat submanifold of a Frobenius manifold $M$ with
$$e\vert_N \in TN;$$
$$TN \circ TN \subseteq TN;$$
$$\mathcal{E}\vert_N \in TN.$$
Then $N$ is a natural Frobenius submanifold
\end{theorem}

Neither $e\vert_N \in TN$ nor $\mathcal{E}\vert_N \in TN$ is
necessary, we will construct examples of natural Frobenius
submanifolds such that $e\vert_N \neq e_N$ and $\mathcal{E}\vert_N
\neq \mathcal{E}_N$.

\begin{example}\label{example1}
Let $(N,g,\ast,e_{N}, \mathcal{E}_{N})$ to be a Frobenius manifold
of dimension~$n$ with constant $D\neq0$, and let $\mathcal{A}$ be
the affine line. Define a new Frobenius manifold $M=N\times
\mathcal{A}$ as follows:

Let $z$ be the coordinate of $\mathcal{A}$ and choose a metric
$\eta$ on $\mathcal{A}$ such that $\eta(\partial_
z,\partial_z)=1$. We define a new metric $\overline g$ on $M$ to
be the direct sum of $g$ and $\eta$. Let $\overline{\nabla}$ be
the Levi-Civita connection of $\overline g$. Then
$\overline{\nabla}$ is just the direct sum of $\nabla$ and $d$,
where $d$ is the Levi-Civita connection of $\eta$. Now define a
product $\circ$: For any $X, Y\in TN$,
\begin{align*}
X\circ Y:&= X\ast Y;\\
X\circ \partial _z&=0;\\
\partial _z\circ \partial _z&=\frac{2}{D}\partial _z.
\end{align*}
Finally we define the unit element $e$ and the Euler vector field $\mathcal{E}$:
\begin{align*}
e&=e_N+\frac{D}{2}\partial_z;\\
\mathcal{E}&=\mathcal{E}_N+\frac{D}{2}z\partial_z.
\end{align*}
It is easy to see $(M,\overline{g},\circ,e, \mathcal{E})$ is a
Frobenius manifold.  Now we embed $N$ to $M$:
$$\iota:N\longrightarrow M, \quad P\longmapsto (P, 1).$$
Then we get a natural Frobenius submanifold $N\times \lbrace
1\rbrace$ with $e_N^\bot\neq 0$ and ${\mathcal{E}_N^\bot}\neq 0$.
\end{example}

\subsection{Aim of the paper}
\label{subsection0c} The paper will give a necessary and
sufficient condition for a submanifold to be a natural Frobenius
manifold and classify the natural Frobenius hypersurfaces.

Let us first recall a known result in differential geometry which
will explain the notation in results below:

\begin{theorem}[\cite{Will}]\label{theorem5}
Let $M$ be a manifold with a metric $\overline g$,
$\overline{\nabla}$ is the Levi-Civita connection of $\overline
g$, and let $N$ be an arbitrary submanifold such that the
restricted metric $g$ is non-degenerate. Then for all $W,X,Y,Z\in
TN$ and normal vectors $\xi,\eta \in TN^\bot$, w.r.t. the
decomposition $TN \oplus TN^{\bot}$, we have:

Gauss formula:
\begin{equation*}
\overline{\nabla}_X Y=\nabla_X Y+h(X,Y),
\end{equation*}

Weingarten formula:
\begin{equation*}
\overline{\nabla}_X \xi =-A_\xi X+\nabla^\bot_X \xi,
\end{equation*}
where $\nabla_X Y:=pr(\overline{\nabla}_X Y)$, $-A_\xi X:=
pr(\overline{\nabla}_X \xi)$.

Here $h$ is called the second fundamental form and $A$ is called
the shape operator, which are related by
\begin{equation*}
\<h(X,Y),\xi \>=\<A_\xi X,Y\>,
\end{equation*}
for any $X, Y \in TN$, $\xi \in TN^\bot$.
\end{theorem}

We have the following result for any submanifold.

\begin{theorem}\label{theorem A}
Let $(M,\overline{g},\circ,e, \mathcal{E})$ be a Frobenius
manifold. Let $N$ be a submanifold of $M$ such that $g:={\overline
g}\mid_{TN}$ is non-degenerate and flat. Then the following
properties are equivalent:

(1) N is a natural Frobenius submanifold of $M$;

(2) \begin{align}
\nabla {e_N}&=0;\label{flateN}\\
TN \circ TN &\subseteq TN;\label{closecirc}\\
\exists \lambda \in \mathbb{C},\quad
A_{\mathcal{E}_N^\bot} &=\lambda \cdot \Id,\label{shapeE}
\end{align}
\end{theorem}

We classify the natural Frobenius hypersurfaces in the following:

\begin{proposition}\label{proposition B}
With the assumptions of Theorem \ref{theorem A}, assume moreover
that $N$ is a hypersurface.

(a) If $e$ is not tangent to $N$, then the following are equivalent:

\begin{enumerate}
\item
$N$ is natural Frobenius submanifold of $M$;

\item
$TN\circ TN \subseteq TN, \overline \nabla=\nabla$.
\end{enumerate}

(b) If $e$ is tangent to $N$, then the following are equivalent:
\begin{enumerate}
\item
$N$ is natural Frobenius submanifold of $M$;
\item
$TN\circ TN \subseteq TN$ and either $ \overline \nabla=\nabla$ or $\mathcal{E}$ is tangent to $N$.
\end{enumerate}
\end{proposition}

We will give some examples in the last section to show that this
classification can not be generalized to all submanifolds.

\subsubsection*{Acknowledgement} The author would like to thank ENS,
Paris for its hospitality during the academic year of 2006--2007.
The author also thanks Jianxun Hu for pointing me to the work of
Strachan. Last but most, the author would like to thank professor
Sabbah for his patient helps and valuable suggestions.

\section{General dimension}\label{section1}

\renewcommand{\theequation}
{1.\arabic{equation}} \setcounter{equation}{0}

In this section we mainly give a necessary and sufficient
condition for a submanifold to be a natural Frobenius manifold.
Because of the equivalence between Frobenius structure and Saito
structure with metric, we will see the sufficient condition from
these two equivalent point of view.

\begin{lemma}\label{lemma C} Let $(M,\overline{g},\circ,e, \mathcal{E})$ be a Frobenius
manifold. Let $N$ be a submanifold of $M$ such that $TN \circ TN
\subseteq TN$, then $TN \circ TN^{\bot} \subseteq TN^{\bot}$.
\end{lemma}
\begin{proof}[The proof of Lemma \ref{lemma C}]Because $(M,\overline{g},\circ,e,
\mathcal{E})$ is a Frobenius manifold, so we have the relation
$$\<U \circ V,W\> = \<U,V \circ W\>,$$
for all $U,V,W \in TM$. If $TN \circ TN \subseteq TN$, then
$\forall U,V \in TN, \zeta \in TN^{\bot}$ we have
$$\<U \circ \zeta,V\> = \<\zeta,U \circ V\> = 0.$$
So $TN \circ TN^{\bot} \subseteq TN^{\bot}$.
\end{proof}

\begin{proof}[Proof of Theorem A]
(1) $\Longrightarrow$ (2). Because $N$ is a Frobenius submanifold
of $M$, there exist two constants $D$ and $D_N$, such that:
\begin{equation*}
\overline{\nabla}\mathcal{E}+(\overline{\nabla}\mathcal{E})^*=D \cdot \Id;
\end{equation*}
\begin{equation*}
\nabla\mathcal{E}_N+(\nabla\mathcal{E}_N)^*=D_N \cdot \Id.
\end{equation*}
For any $U,V\in TN$, we have:
\begin{equation*}
\< \overline{\nabla}_U \mathcal{E},V\>+\<\overline{\nabla}_V
\mathcal{E},U\>=D\cdot \<U,V\>
\end{equation*}
Computing the left hand side of the above equality we get:
\begin{eqnarray*}
\<\overline{\nabla}_U \mathcal{E},V\>+\<\overline{\nabla}_V
\mathcal{E},U\> & = & \<\overline{\nabla}_U
\mathcal{E}_N+\overline{\nabla}_U \mathcal{E}_N^ \bot,V\>
+\<\overline{\nabla}_V \mathcal{E}_N+\overline{\nabla}_V \mathcal{E}_N^ \bot,U\> \\
& = & \<{\nabla}_U \mathcal{E}_N + {\nabla}^\bot _U \mathcal{E}_N^ \bot+h(U,\mathcal{E}_N) - A_{\mathcal{E}_N^ \bot}U,V\> \\
& & +
\<{\nabla}_V \mathcal{E}_N+{\nabla}^\bot _V \mathcal{E}_N^ \bot + h(V,\mathcal{E}_N)-A_{\mathcal{E}_N^ \bot}V,U\> \\
& = & \<{\nabla}_U \mathcal{E}_N-A_{\mathcal{E}_N^ \bot}U,V\>+
\<{\nabla}_V \mathcal{E}_N - A_{\mathcal{E}_N^ \bot}V,U\> \\
& = & \<{\nabla}_U \mathcal{E}_N,V\>+ \<{\nabla}_V \mathcal{E}_N,
U\>
-\<A_{\mathcal{E}_N^ \bot}U,V\> -
\<A_{\mathcal{E}_N^ \bot}V,U\> \\
& = & \<{\nabla}_U \mathcal{E}_N,V\>+ \<{\nabla}_V \mathcal{E}_N,
U\>
-2\<h(U,V),{\mathcal{E}_N^\bot}\> \\
& = & D_N \cdot \<U,V\>-2\<h(U,V),{\mathcal{E}_N^\bot}\> \\
& = & D \cdot \<U,V\>
\end{eqnarray*}
So we get:
\begin{equation*}
\<h(U,V),{\mathcal{E}_N^\bot}\>=\frac{D_N-D}{2}\cdot \<U,V\>,
\end{equation*}
i.e.,
\begin{equation*}
\<A_{\mathcal{E}_N^\bot}U,V\>=\frac{D_N-D}{2}\cdot \<U,V\>.
\end{equation*}
So
\begin{equation*}
A_{\mathcal{E}_N^\bot}=\frac{D_N-D}{2}\cdot \Id
\end{equation*}

The other two equalities $\nabla {e_N}=0$ and $TN \circ TN
\subseteq TN$ hold because $N$ is the natural Frobenius
submanifold of $M$.

(2) $\Longrightarrow$ (1)

We will give two methods to prove the sufficient condition. In the
first method, we use the flat holomorphic local coordinates to
prove that Condition (2) induces a Saito structure with metric on
$M$. The most difficult part in this method is the flatness of the
structure connection $\widetilde{\nabla}$. The second method is
more global. We prove that Condition (2) induces a Frobenius
structure on $M$. In this method every thing is more obvious
except the relation $\mathcal{L}_{\mathcal{E}_N}( \circ )= \circ$.

\subsubsection*{First method: Saito structure}

Suppose $TN \circ TN \subseteq TN$ and there exists a constant $\lambda \in \mathbb{C}$ such that $A_{\mathcal{E}_N^\bot} =\lambda \cdot \Id$. So the restricted Higgs field $\Phi|_{TN}$ is a Higgs field on $TN$, where $\Phi$ is defined by $\Phi _X Y:= -X\circ Y$,for any $X, Y \in TM$.

The structure connections $\widetilde{\overline{\nabla}}$ on
$M\times \mathbf{P}^1$ and $\widetilde{\nabla}$ on $N\times
\mathbf{P}^1$ are defined by
\begin{equation*}
\widetilde{\overline{\nabla}}:=\pi^*\overline \nabla+ \frac{\pi^*\Phi}{z}-(\frac{\Phi(\mathcal{E})}{z}+\overline \nabla \mathcal{E}) \frac{dz}{z};
\end{equation*}
\begin{equation*}
\widetilde{\nabla}:=\pi^* \nabla+ \frac{\pi^*(\Phi|_{TN})}{z}-(\frac{\Phi(\mathcal{E}_N)}{z}+ \nabla \mathcal{E}_N) \frac{dz}{z};
\end{equation*}
We will show that the induced structure $(\nabla, \Phi|_{TN}, e_N,
\mathcal{E}_N, g)$ on $N$ is a Saito structure with metric.

\noindent $S1)$ \emph{Existence of flat unit field}. Because of
$TN \circ TN \subseteq TN$. From Lemma \ref{lemma C}, we know $TN
\circ TN^\bot \subseteq TN^\bot$, so for any $U \in TN$, we have:
$$U=U \circ e=U \circ e_N + U \circ e_N^\bot.$$
So
$$U \circ e_N^\bot = U-U \circ e_N \in TN\cap TN^\bot=\lbrace 0 \rbrace.$$
So
$$U=U \circ e_N,$$ for any $U \in TN$.
i.e., ${\Phi_{e_N}}{|_{TN}}=-\Id$.
$\nabla {e_N}=0$ show that the unit vector field $e_N$ is $\nabla$-flat.

\noindent $S2)$ \emph{flatness of the structure connection
$\widetilde{\nabla}$.} Denote by
$\widetilde{\overline{\mathcal{R}}}$ the curvature of
$\widetilde{\overline{\nabla}}$ and by $\widetilde{\mathcal{R}}$
the curvature of $\widetilde{\nabla}$. Because $M$ is Frobenius
manifold, $\widetilde{\overline {\mathcal{R}}}=0$. For any $U,V,
W \in TM$, we have:
$$\widetilde{\overline {\mathcal{R}}}(U,V)W=0.$$
Computing the left hand side of the above equality:
\begin{eqnarray*}
\widetilde{\overline {\mathcal{R}}}(U,V)W
& = & \overline{R}(U,V)W \\
& & + \frac{1}{z} \lbrace U\circ \overline \nabla_V W - V \circ \overline \nabla_U W - [U,V]\circ W + \overline \nabla_U (V\circ W) - \overline \nabla_V (U\circ W) \rbrace \\
& & + \frac{1}{z^2} \lbrace U \circ (V\circ W) - V \circ (U \circ W) \rbrace,
\end{eqnarray*}
where $\overline{R}$ is the curvature of $\overline \nabla$. So we get:
\begin{equation}
U\circ \overline \nabla_V W - V \circ \overline \nabla_U W -
[U,V]\circ W + \overline \nabla_U (V\circ W) - \overline \nabla_V
(U\circ W)=0.
\end{equation}
because of Lemma \ref{lemma C} we have, for $\forall U,
V,W \in TN$
\begin{eqnarray*}
& & \pr \lbrace U\circ
\overline \nabla_V W - V \circ \overline \nabla_U W - [U,V]\circ W
+ \overline \nabla_U (V\circ W) - \overline \nabla_V (U\circ
W)\rbrace \\
& = & U\circ \nabla_V W - V \circ \nabla_U W - [U,V]\circ W + \nabla_U (V\circ W) - \nabla_V (U\circ W) \\
& = & 0.
\end{eqnarray*}
However,
\begin{eqnarray*}
\widetilde{ {\mathcal{R}}}(U,V)W
& = & R(U,V)W \\
& & + \frac{1}{z} \lbrace U\circ \nabla_V W - V \circ \nabla_U W -
[U,V]\circ W + \nabla_U (V\circ W) - \nabla_V (U\circ W) \rbrace
\\
& & + \frac{1}{z^2} \lbrace U \circ (V\circ W) - V \circ (U \circ
W) \rbrace,
\end{eqnarray*}
where $R$ is the curvature of $\nabla$. So
\begin{equation*}
\widetilde{\mathcal{R}}(U,V)W=0,
\end{equation*}
for any $U,V,W \in TN$

Now the only other equality to be checked is
$\widetilde{{\mathcal{R}}}(z\frac{d}{dz},U) V=0$. Calculating
directly, we get:
\begin{equation*}
\widetilde{{\mathcal{R}}}(z\frac{d}{dz},U) V =- \nabla_U \nabla_V {\mathcal{E}_N}+ \nabla_{\nabla_U V} \mathcal{E}_N.
\end{equation*}
for any $U,V \in TN$.

Suppose $t^1,t^2,\dots,t^m$ is the flat coordinate of $(M, \overline{\nabla})$, $\tau^1, \tau^2,\dots, \tau^n$ is the flat coordinate of $N$, $\widetilde{\mathcal{R}}$ is a tensor, so we just check it for base elements $\partial_\alpha$. So we just need to check:
\begin{equation*}
\partial_{\tau^\alpha}\partial_{\tau^\beta}\mathcal{E}_{N}^ \gamma=0,
\end{equation*}
for all $\alpha,\beta,\gamma \in \lbrace1,2,\dots,n\rbrace$.

Locally, $\mathcal{E}=\mathcal{E}^i\partial_{t^i},
\mathcal{E}_{N}=\mathcal{E}_{N}^\alpha \partial_{\tau^\alpha}$.
Choose the local frame of $TN^\bot$, denoted by
$\partial_{\nu^{\tilde\alpha}}$, such that
\begin{equation*}
\partial_{t^i}=A_i^\alpha\partial_{\tau^\alpha}+n^{\tilde\alpha}_i\partial_{\nu^{\tilde\alpha}}.
\end{equation*}
and
\begin{equation*}
\<\partial_{\nu^{\tilde\alpha}}, \partial_{\nu^{\tilde\beta}}\> =
\eta_{{\tilde\alpha}{\tilde\beta}}
\end{equation*}
where $\eta_{{\tilde\alpha}{\tilde\beta}}$ are constant with
$\eta_{{\tilde\alpha}{\tilde\beta}}=\epsilon({\tilde\alpha})
\delta_{{\tilde\alpha}{\tilde\beta}}$ with $\epsilon({\tilde\alpha})=\pm 1$.

Using the metrics $\overline g$ and $g$ we get:
\begin{equation*}
A_i^\alpha=\overline g_{ij} g^{\alpha \beta} \frac{\partial_{t^j}}{\partial_{\tau^ \beta}}
\end{equation*}
From the definition of $\mathcal{E}_{N}$, we get $\mathcal{E}_{N}^\alpha = \mathcal{E}^i\mid_NA_i^\alpha$.

Computing $\frac {\partial \mathcal{E}_{N}^\gamma}{\partial_{\tau^\beta}}$ directly we get:
\begin{equation*}
\frac {\partial \mathcal{E}_{N}^\gamma}{\partial_{\tau^\beta}}=\frac {\partial \mathcal{E}^i }{ \partial \tau_ \beta}\overline{g}_{ij}g^{\gamma\delta} \frac{\partial t^j }{ \partial \tau_ \delta}+\mathcal{E}^i\overline{g}_{ij}g^{\gamma\delta}\frac {\partial^2t^j}{\partial_{\tau^\beta}\partial_{\tau^\delta}}
\end{equation*}
On the other hand,
\begin{eqnarray*}
\mathcal{E}^i\overline{g}_{ij}g^{\gamma\delta}\frac
{\partial^2t^j}{\partial_{\tau^\beta}\partial_{\tau^\delta}}
& = & \mathcal{E}^i \<\partial_{t^i},\partial_{t^j}\>g^{\gamma\delta}\frac {\partial^2t^j}{\partial_{\tau^\beta}\partial_{\tau^\delta}} \\
& = & \big \<\mathcal{E}^i \partial_{t^i}, \frac {\partial^2t^j}{\partial_{\tau^\beta}\partial_{\tau^\delta}}\partial_{t^j} \big \>g^{\gamma\delta} \\
& = & g^{\gamma\delta} \<\mathcal{E}, \overline \nabla_{\partial_{\tau^\beta}} \partial_{\tau^\delta} \> \\
& = & g^{\gamma\delta} \<\mathcal{E},\nabla_{\partial_{\tau^\beta}} \partial_{\tau^\delta} + h(\partial_{\tau^\beta},\partial_{\tau^\delta})\>
\end{eqnarray*}
But $\tau^1, \tau^2,\dots, \tau^n$ is the flat coordinate of $N$, so
\begin{eqnarray*}
\mathcal{E}^i\overline{g}_{ij}g^{\gamma\delta}\frac
{\partial^2t^j}{\partial_{\tau^\beta}\partial_{\tau^\delta}}
& = & g^{\gamma\delta}\<\mathcal{E},h(\partial_{\tau^\beta},\partial_{\tau^\delta})\>\\
& = & g^{\gamma\delta}\<\mathcal{E}_N^\bot,h(\partial_{\tau^\beta},\partial_{\tau^\delta})\>\\
& = & g^{\gamma\delta}\<A_{\mathcal{E}_N^\bot}\partial_{\tau^\beta},\partial_{\tau^\delta}\>\\
& = & \lambda g^{\gamma \delta}g_{\delta \beta}.
\end{eqnarray*}
The last equality holds because $A_{\mathcal{E}_N^\bot} =\lambda \cdot \Id$.

So we have:
\begin{equation*}
\frac {\partial \mathcal{E}_{N}^\gamma}{\partial_{\tau^\beta}}=\frac {\partial \mathcal{E}^i }{ \partial \tau_ \beta}\overline{g}_{ij}g^{\gamma\delta} \frac{\partial t^j }{ \partial \tau_ \delta}+\lambda g^{\gamma \delta}g_{\delta \beta}
\end{equation*}
Similarly computing we get:
\begin{equation*}
\frac{\partial^2\mathcal{E}_{N}^\gamma}{\partial_{\tau^\alpha}\partial_{\tau^\beta}}=
\frac {\partial^2\mathcal{E}^i}{\partial_{\tau^\alpha}\partial_{\tau^\beta}} \overline g_{ij} g^{\gamma \delta} \frac{\partial_{t^j}}{\partial_{\tau^ \delta}}
+ \frac {\partial^2 t^j}{\partial_{\tau^\alpha}\partial_{\tau^\delta}} \overline g_{ij} g^{\gamma \delta} \frac{\partial_{\mathcal{E}^i}}{\partial_{\tau^ \delta}}
\end{equation*}
We will prove that the right hand side of the above equality
vanishes.

Calculating the first term of the right hand side we get:
\begin{eqnarray*}
\frac
{\partial^2\mathcal{E}^i}{\partial_{\tau^\alpha}\partial_{\tau^\beta}}
\overline g_{ij} g^{\gamma \delta}
\frac{\partial_{t^j}}{\partial_{\tau^ \delta}}
& = & \frac {\partial^2\mathcal{E}^i}{\partial_{\tau^\alpha}\partial_{\tau^\beta}} \<\partial_{t^i}, \partial_{t^j}\> g^{\gamma \delta} \frac{\partial_{t^j}}{\partial_{\tau^ \delta}}\\
& = & \big \< \frac
{\partial^2\mathcal{E}^i}{\partial_{\tau^\alpha}\partial_{\tau^\beta}}
\partial_{t^i}, \frac{\partial t^j}{\partial_{\tau^ \delta}}
\partial_{t^j}\big \> g^{\gamma \delta}.
\end{eqnarray*}
Claim:
$\overline \nabla_{h(\partial_{\tau^\alpha},\partial_{\tau^\beta})}\mathcal{E}$
$=\frac{\partial^2\mathcal{E}^i}{\partial_{\tau^\alpha}\partial_{\tau^\beta}} \partial_{t^i}$

In fact, because $M$ is a Frobenius manifold, we have $\widetilde{\overline {\mathcal{R}}}=0$,
so we have
\begin{eqnarray*}
\widetilde{\overline {\mathcal{R}}}(z\frac{d}{dz},U) V =- \overline \nabla_U \overline \nabla_V {\mathcal{E}}+ \overline \nabla_{\overline \nabla_U V} \mathcal{E}=0,
\end{eqnarray*}
for all $U,V \in TN$. By this equality we get:
\begin{eqnarray*}
\overline\nabla_{h(\partial_{\tau^\alpha},\partial_{\tau^\beta})}\mathcal{E}
& = & \overline \nabla _{\overline \nabla_{\partial_{\tau^\alpha}} \partial_{\tau^\beta}} \mathcal{E}\\
& = & \overline \nabla_{\partial_{\tau^\alpha}} \overline \nabla_{\partial_{\tau^\beta}} \mathcal{E} \\
& = & \frac{\partial^2\mathcal{E}^i}{\partial_{\tau^\alpha}\partial_{\tau^\beta}} \partial_{t^i}.
\end{eqnarray*}
So
\begin{equation*}
\frac
{\partial^2\mathcal{E}^i}{\partial_{\tau^\alpha}\partial_{\tau^\beta}}
\overline g_{ij} g^{\gamma
\delta}\frac{\partial_{t^j}}{\partial_{\tau^ \delta}}
=\frac{\partial^2\mathcal{E}_{N}^\gamma}{\partial_{\tau^\alpha}\partial_{\tau^\beta}}=\<\overline\nabla_{h(\partial_{\tau^\alpha},\partial_{\tau^\beta})}\mathcal{E},\partial_{\tau^\delta}\>g^{\gamma
\delta}
\end{equation*}
Similarly for the second term:
\begin{equation*}
\frac {\partial^2 t^j}{\partial_{\tau^\alpha}\partial_{\tau^\delta}} \overline g_{ij} g^{\gamma \delta} \frac{\partial {\mathcal{E}^i}}{\partial_{\tau^ \delta}}
=\<\overline\nabla_{\partial_{\tau^\beta}}\mathcal{E},h(\partial_{\tau^\delta},\partial_{\tau^\alpha})\>g^{\gamma\delta}
\end{equation*}
We simplify the equality to be:
\begin{equation*}
\frac{\partial^2\mathcal{E}_{N}^\gamma}{\partial_{\tau^\alpha}\partial_{\tau^\beta}}=\<\overline\nabla_{h(\partial_{\tau^\alpha},\partial_{\tau^\beta})}\mathcal{E},\partial_{\tau^\delta}\>g^{\gamma
\delta}+\<\overline\nabla_{\partial_{\tau^\beta}}\mathcal{E},h(\partial_{\tau^\delta},\partial_{\tau^\alpha})\>g^{\gamma\delta}
\end{equation*}

That $M$ is a Frobenius manifold also implies that there exists a
constant $D$ such that:
\begin{equation*}
\overline{\nabla}\mathcal{E}+(\overline{\nabla}\mathcal{E})^*=D\cdot \Id;
\end{equation*}
So
\begin{eqnarray*}
\<\overline\nabla_{h(\partial_{\tau^\alpha},\partial_{\tau^\beta})}\mathcal{E},\partial_{\tau^\delta}\>
& = & D
\cdot\<h(\partial_{\tau^\alpha},\partial_{\tau^\beta}),\partial_{\tau^\delta}\>
-\<h(\partial_{\tau^\alpha},\partial_{\tau^\beta}),\overline\nabla_{\partial_{\tau^\delta}}\mathcal{E}\>\\
& = & -\<h(\partial_{\tau^\alpha},\partial_{\tau^\beta}),\overline\nabla_{\partial_{\tau^\delta}}\mathcal{E}\>.
\end{eqnarray*}
Because $\overline \nabla \overline g=0$, we get:
\begin{eqnarray*}
-\<h(\partial_{\tau^\alpha},\partial_{\tau^\beta}),
\overline\nabla_{\partial_{\tau^\delta}}\mathcal{E}\>
& = &
-\partial_{\tau^\delta}(\<h(\partial_{\tau^\alpha},\partial_{\tau^\beta}),
\mathcal{E}\>)
+ \<\overline \nabla_{\partial_{\tau^\delta}}h(\partial_{\tau^\alpha},\partial_{\tau^\beta}),\mathcal{E}\>\\
& = & -\partial_{\tau^\delta}(\<h(\partial_{\tau^\alpha},\partial_{\tau^\beta}), \mathcal{E}_N^\bot\>)
+\<\overline \nabla_{\partial_{\tau^\delta}}h(\partial_{\tau^\alpha},\partial_{\tau^\beta}),\mathcal{E}\>\\
& = & -\partial_{\tau^\delta}(\lambda g_{\alpha \beta})
+\<\overline \nabla_{\partial_{\tau^\delta}}h(\partial_{\tau^\alpha},\partial_{\tau^\beta}),\mathcal{E}\>\\
& = & \<\overline \nabla_{\partial_{\tau^\delta}}h(\partial_{\tau^\alpha},\partial_{\tau^\beta}),\mathcal{E}\>.
\end{eqnarray*}
Similarly we get:
\begin{eqnarray*}
\<\overline\nabla_{\partial_{\tau^\beta}}\mathcal{E},h(\partial_{\tau^\delta},\partial_{\tau^\alpha})\>
=-\<\overline \nabla_{\partial_{\tau^\beta}}h(\partial_{\tau^\alpha},\partial_{\tau^\delta}),\mathcal{E}\>
\end{eqnarray*}
Then we have the equality:
\begin{equation*}
\frac{\partial^2\mathcal{E}_{N}^\gamma}{\partial_{\tau^\alpha}\partial_{\tau^\beta}}
=g^{\delta\gamma}\<\mathcal{E},\overline\nabla_{\partial_{\tau^\delta}}h(\partial_{\tau^\alpha},\partial_{\tau^\beta})-\overline\nabla_{\partial_{\tau^\beta}}h(\partial_{\tau^\delta},\partial_{\tau^\alpha})\>
\end{equation*}
However the right hand side of this equality vanishes because $\overline\nabla_{\partial_{\tau^\delta}}h(\partial_{\tau^\alpha},\partial_{\tau^\beta})$ is totally symmetric in $\alpha,\beta,\delta:$
\begin{eqnarray*}
\overline\nabla_{\partial_{\tau^\delta}}h(\partial_{\tau^\alpha},\partial_{\tau^\beta})
& = & \overline\nabla_{\partial_{\tau^\delta}}\overline\nabla_{\partial_{\tau^\beta}}\partial_{\tau^\alpha}\\
& = & \overline\nabla_{\partial_{\tau^\beta}}\overline\nabla_{\partial_{\tau^\delta}}\partial_{\tau^\alpha}
+\overline\nabla_{[\partial_{\tau^\delta},\partial_{\tau^\beta}]}\partial_{\tau^\alpha}\\
& = & \overline\nabla_{\partial_{\tau^\beta}}h(\partial_{\tau^\delta},\partial_{\tau^\alpha}).
\end{eqnarray*}
Then we get
$\partial_{\tau^\alpha}\partial_{\tau^\beta}\mathcal{E}_{N}^\gamma=0,
\forall \alpha,\beta,\gamma \in \lbrace1,2,\dots,n\rbrace$. So
$\widetilde{\mathcal{R}}=0$, i.e., the structure connection
$\widetilde{\nabla}$ is integrable.

From $S1)$ and $S2)$ we get $(\nabla, \Phi_{|_{TN}}, e_N, \mathcal{E}_N)$ is a Saito structure (without metric) on $N$.

\noindent $S3)$ \emph{Saito structure $(\nabla, \Phi_{|_{TN}},
e_N, \mathcal{E}_N)$ with metric $g$}. Because $\nabla$ is the
Levi-Civita connection of $g$, so we have:
$$\nabla g = 0.$$
The induced Higgs field $\Phi|_{TN}$ satisfies $\Phi|_{TN} =
(\Phi|_{TN})^*$ w.r.t. $g$ because $\Phi = (\Phi)^*$ w.r.t.
$\overline g$. So we just need to check:

$\exists D_N \in \mathbb{C}$, such that
$$\nabla\mathcal{E}_N+(\nabla \mathcal{E}_N)^*=D_N \cdot \Id.$$
Because $M$ is Frobenius manifold, there exists a constant $D$ such that:
$$\overline \nabla \mathcal{E}+(\overline \nabla \mathcal{E})^*=D \cdot \Id.$$

Computing the left hand side of the above relation as in the proof
of $(1) \Rightarrow (2)$, we get for any $U,V \in TN$
\begin{eqnarray*}
\<\overline{\nabla}_U \mathcal{E},V\>+\<\overline{\nabla}_V
\mathcal{E},U\> & = & \<{\nabla}_U \mathcal{E}_N,V\>+
\<{\nabla}_V \mathcal{E}_N,U\>
-\<A_{\mathcal{E}_N^ \bot}U,V\> \\
& & -
\<A_{\mathcal{E}_N^ \bot}V,U\> \\
& = & \<{\nabla}_U \mathcal{E}_N,V\>+
\<{\nabla}_V \mathcal{E}_N,U\> - 2\lambda \cdot \<U,V\> \\
& = & D\cdot \<U,V\>
\end{eqnarray*}
so
$$\<{\nabla}_U \mathcal{E}_N,V\>+\<{\nabla}_V \mathcal{E}_N,U\>= D\cdot \<U,V\>+2\lambda \cdot \<U,V\>$$
That is to say:
$$\nabla \mathcal{E}_N+(\nabla \mathcal{E}_N)^*=(D+2\lambda)\cdot \Id$$
Take $D_N=D+2\lambda$, we get the
equality:$$\nabla\mathcal{E}_N+(\nabla \mathcal{E}_N)^*=D_N \cdot
\Id.$$

From $S1),S2),S3)$ we know that $(\nabla, \Phi_{|_{TN}}, e_N, \mathcal{E}_N, g)$ is a Saito structure on $N$.

\subsubsection*{Second method: Frobenius manifold structure}

Consider the quadruple $(N, \circ, g, e_N, \mathcal{E}_N)$.

\noindent $F1)$ From the assumption we know that $g$ is flat. Just
like the proof of $S1)$, we get the unit vector field $e_N$ is
$\nabla$-flat.

\noindent $F2)$ In other hand, $TN \circ TN \subseteq TN$ also implies that
\begin{equation*}
\< U \circ V,W\> = \< U,V\circ W \>, \forall U,V,W \in TN.
\end{equation*}

\noindent $F3)$ Now define a new 3-tensor
\begin{equation*}
c_{N}(U,V,W):= \<U \circ V,W\>.
\end{equation*}
It is easy to see that $c_{N}$ is the restricted tensor of $c$ to
$TN\otimes TN \otimes TN$, where
\begin{equation*}
c(\overline U, \overline V, \overline W):= \< \overline U \circ
\overline V, \overline W\>, \forall \overline U, \overline V.
\overline W \in TM.
\end{equation*}
$M$ is a Frobenius manifold, so the 4-tensor $(\overline
\nabla_{\overline W^{'}} c)(\overline U, \overline V, \overline
W)$ is symmetric in the four vector fields $\overline U, \overline
V, \overline W, \overline W^{'} \in TM$.

So for any $U,V,W,W^{'} \in TN$ we have
\begin{align*}
(\overline \nabla_{W^{'}} c)&(U,V,W)\\
& = W^{'} (c(U,V,W))- c(\overline \nabla_{W^{'}} U,V,W) - c(U, \overline \nabla_{W^{'}} V,W)- c(U, v, \overline \nabla_{W^{'}} W) \\
& = W^{'} (c_{N} (U,V,W)) - c(\nabla_{W^{'}} U,V,W)- c(U, \nabla_{W^{'}} V,W)- c(U, v, \nabla_{W^{'}} W) \\
& - c(h(W^{'},U),V,W)- c(U, h(W^{'},V),W)- c(U, v, h(W^{'},W)) \\
\end{align*}
However for any $U,V,W,W^{'} \in TN$ we have
\begin{equation*}
c(h(W^{'},U),V,W)= \<h(W^{'},U) \circ V,W\> = \<h(W^{'},U),W \circ V \> = 0.
\end{equation*}

So for any $U,V,W,W^{'} \in TN$ we get
\begin{align*}
(\overline \nabla_{W^{'}}& c)(U,V,W)\\
& = W^{'} (c_{N} (U,V,W)) - c(\nabla_{W^{'}} U,V,W)- c(U, \nabla_{W^{'}} V,W)- c(U, v, \nabla_{W^{'}} W) \\
& = W^{'} (c_{N} (U,V,W)) - c_{N}(\nabla_{W^{'}} U,V,W)- c_{N}(U, \nabla_{W^{'}} V,W)- c_{N}(U, v, \nabla_{W^{'}} W) \\
& = (\nabla_{W^{'}} c_{N})(U,V,W).
\end{align*}
But we know that the 4-tensor $(\overline \nabla_{\overline W^{'}}
c)(\overline U, \overline V, \overline W)$ is symmetric in the
four vector fields $\overline U, \overline V,\overline W,
\overline W^{'} \in TM$. Specially, it is symmetric in the four
vector fields $U,V,W,W^{'} \in TN$, i.e., we get
\begin{equation*}
(\nabla_{W^{'}} c_{N})(U,V,W)
\end{equation*}
is symmetric in the four vector fields $U,V,W,W^{'} \in TN$.

\noindent $F4)$ Now consider the vector field $\mathcal{E}_N$

$M$ is a Frobenius manifold, so there exists a constant $D$ such that
\begin{equation*}
\overline{\nabla}\mathcal{E}+(\overline{\nabla}\mathcal{E})^*=D \cdot \Id.
\end{equation*}
Computing the left hand side of the above relation as in
$(1)\Rightarrow(2)$, together with the condition
$A_{\mathcal{E}_N^\bot} =\lambda \cdot \Id$, we have for any $U,V \in TN$
\begin{eqnarray*}
\<\overline{\nabla}_U \mathcal{E},V\>+\<\overline{\nabla}_V
\mathcal{E},U\> & = & \<{\nabla}_U \mathcal{E}_N,V\>+
\<{\nabla}_V \mathcal{E}_N,U\>
-\<A_{\mathcal{E}_N^ \bot}U,V\>
-\<A_{\mathcal{E}_N^ \bot}V,U\> \\
& = & \<{\nabla}_U \mathcal{E}_N,V\>+
\<{\nabla}_V \mathcal{E}_N,U\> - 2\lambda \cdot \<U,V\> \\
& = & D\cdot \<U,V\>,
\end{eqnarray*}
so
$$\<{\nabla}_U \mathcal{E}_N,V\>+\<{\nabla}_V \mathcal{E}_N,U\>= D\cdot \<U,V\>+2\lambda \cdot \<U,V\>$$
That is to say:
$$\nabla \mathcal{E}_N+(\nabla \mathcal{E}_N)^*=(D+2\lambda)\cdot \Id$$
Take $D_N=D+2\lambda$, we get the
equality:$$\nabla\mathcal{E}_N+(\nabla \mathcal{E}_N)^*=D_N \cdot
\Id.$$

Because $g$ is flat, this relation is equivalent to
\begin{equation*}
\mathcal{L}_{\mathcal{E}_N} (g) = D_N \cdot g.
\end{equation*}

\noindent $F5)$ We will prove
\begin{equation*}
\mathcal{L}_{\mathcal{E}_N} (\circ ) = \circ.
\end{equation*}
Modulo the relation $\nabla (\Phi |_{TN})=0$ this is equivalent to the relation:
\begin{equation*}
\nabla_U(V \circ \mathcal{E}_N)- (\nabla_U V )\circ \mathcal{E}_N + U \circ \nabla_V \mathcal{E}_N - \nabla_{U \circ V} \mathcal{E}_N = U \circ V,
\end{equation*}
for any $U,V \in TN$.

However $M$ is a Frobenius manifold, so we have
\begin{equation*}
\overline \nabla_U(V \circ \mathcal{E})- (\overline \nabla_U V )\circ \mathcal{E} + U \circ \overline \nabla_V \mathcal{E}- \overline \nabla_{U \circ V} \mathcal{E} = U \circ V,
\end{equation*}
for any $U,V \in TN$.

We compute the l.h.s. of this equality and get
\begin{multline*}
\overline \nabla_U(V \circ \mathcal{E})- (\overline \nabla_U V
)\circ \mathcal{E} + U \circ \overline \nabla_V \mathcal{E}-
\overline \nabla_{U \circ V} \mathcal{E} \\
= \overline \nabla_U(V \circ \mathcal{E}_N)- (\overline \nabla_U V )\circ \mathcal{E}_N + U \circ \overline \nabla_V \mathcal{E}_N- \overline \nabla_{U \circ V} \mathcal{E}_N \\
+ \overline \nabla_U(V \circ \mathcal{E}_N^\bot)- (\overline \nabla_U V )\circ \mathcal{E}_N^\bot + U \circ \overline \nabla_V \mathcal{E}_N^\bot- \overline \nabla_{U \circ V} \mathcal{E}_N^\bot
\end{multline*}
Computing first term
\begin{multline*}
\pr (\overline \nabla_U(V \circ \mathcal{E}_N)- (\overline \nabla_U V )\circ \mathcal{E}_N + U \circ \overline \nabla_V \mathcal{E}_N- \overline \nabla_{U \circ V} \mathcal{E}_N) \\
= \nabla_U(V \circ \mathcal{E}_N)- (\nabla_U V )\circ \mathcal{E}_N + U \circ \nabla_V \mathcal{E}_N - \nabla_{U \circ V} \mathcal{E}_N
\end{multline*}
So we just need to prove
\begin{equation*}
\pr(\overline \nabla_U(V \circ \mathcal{E}_N^\bot)- (\overline
\nabla_U V )\circ \mathcal{E}_N^\bot + U \circ \overline \nabla_V
\mathcal{E}_N^\bot- \overline \nabla_{U \circ V}
\mathcal{E}_N^\bot) =0.
\end{equation*}
computing directly we find
\begin{align*}
\pr(\overline \nabla_U(V \circ \mathcal{E}_N^\bot)&- (\overline
\nabla_U V )\circ \mathcal{E}_N^\bot + U \circ \overline \nabla_V
\mathcal{E}_N^\bot- \overline \nabla_{U \circ V}
\mathcal{E}_N^\bot) \\
& = - A_{V \circ \mathcal{E}_N^\bot}(U) - \pr(h(U,V) \circ \mathcal{E}_N^\bot) - U \circ A_{\mathcal{E}_N^\bot}(V) + A_{\mathcal{E}_N^\bot}(U \circ V) \\
& = - A_{V \circ \mathcal{E}_N^\bot}(U) - \pr( h(U,V) \circ \mathcal{E}_N^\bot) - \lambda U \circ V + \lambda U \circ V \\
& = - A_{V \circ \mathcal{E}_N^\bot}(U) - \pr( h(U,V) \circ
\mathcal{E}_N^\bot ).
\end{align*}
The second equality holds because $A_{\mathcal{E}_N^\bot } = \lambda \cdot \Id$.

\emph{Claim}: $- A_{V \circ \mathcal{E}_N^\bot}(U) - \pr( h(U,V)
\circ \mathcal{E}_N^\bot )=0$ for any $U,V \in TN$.

In fact,

$(1^{\circ})$ the structure connection of $M$ is flat because $M$
is Frobenius manifold. So the relation $(1.2)$ holds:
\begin{equation*}
U\circ \overline \nabla_V W - V \circ \overline \nabla_U W - [U,V]\circ W + \overline \nabla_U (V\circ W) - \overline \nabla_V (U\circ W)=0.
\end{equation*}
so the orthogonal part of this coefficient must be zero, i.e.,
\begin{equation*}
U\circ h(V,W)-V \circ h(U,W)+h(U,V\circ W)-h(V,U\circ W)=0.
\end{equation*}
i.e.,
\begin{equation*}
U\circ h(V,W)-V \circ h(U,W)= h(V,U\circ W)-h(U,V\circ W).
\end{equation*}
then we get:
\begin{equation*}
\< U\circ h(V,W)-V \circ h(U,W), \mathcal{E}_N^\bot \> = \< h(V,
U\circ W)-h(U,V\circ W), \mathcal{E}_N^\bot \>.
\end{equation*}
We simply the r.h.s. of this equality
\begin{eqnarray*}
\< h(V,U\circ W)-h(U,V\circ W), \mathcal{E}_N^\bot \>
& = & \< h(V,U\circ W), \mathcal{E}_N^\bot \> - \< h(U,V\circ W), \mathcal{E}_N^\bot \> \\
& = & \<A_{\mathcal{E}_N^\bot}(V),W \circ U\> - \<A_{\mathcal{E}_N^\bot}(U),W \circ V\> \\
& = & \lambda \<V,W \circ U\> - \lambda \<U,W \circ V\> \\
& = & 0.
\end{eqnarray*}
So
$$\< U\circ h(V,W), \mathcal{E}_N^\bot \> = \< V \circ h(U,W), \mathcal{E}_N^\bot \>.$$

$(2^{\circ} )$ For any $W \in TN$, we have
\begin{eqnarray*}
\< A_{V \circ {\mathcal{E}_N^\bot}}U,W \>
& = & \< \overline \nabla_U ({V \circ {\mathcal{E}_N^\bot}}),W \> \\
& = & U \< V \circ \mathcal{E}_N^\bot,W \> - \< V \circ {\mathcal{E}_N^\bot}, \overline \nabla_U W \> \\
& = & 0- \< V \circ {\mathcal{E}_N^\bot}, h(U,W) \> \\
& = & - \< {\mathcal{E}_N^\bot},V \circ h(U,W) \>
\end{eqnarray*}
the second equality holds because $\overline \nabla$ is the Levi-Civita connection of $\overline g$, the last equality holds because the product is compatible to the metric $\overline g$.

Now We consider the second term
\begin{equation*}
\<h(U,V)\circ {\mathcal{E}_N^\bot},W \> = \<{\mathcal{E}_N^\bot}, h(U,V)\circ W \>.
\end{equation*}
for all $U,V,W \in TN$. But in $(1^{\circ})$ We have proved that
$$\< U\circ h(V,W), \mathcal{E}_N^\bot \> = \< V \circ h(U,W), \mathcal{E}_N^\bot \>.$$
So for any $Z \in TN$ we have
\begin{equation*}
\<h(U,V)\circ {\mathcal{E}_N^\bot},W \> + \< A_{V \circ
{\mathcal{E}_N^\bot}}U,W \> = 0,
\end{equation*}
i.e., for any $U,V \in TN$ we have:
$$\pr(A_{V \circ {\mathcal{E}_N^\bot}}U+h(U,V)\circ {\mathcal{E}_N^\bot})=0.$$

So we get
\begin{equation*}
\mathcal{L}_{\mathcal{E}_N} (\circ ) = \circ
\end{equation*}

\noindent $F6)$ From Remark \ref{remark SABB} applied to $N$, we
deduce that
\begin{equation*}
\nabla(\nabla \mathcal{E}_N)=0.
\end{equation*}
So $(N, \circ, g, e_N, \mathcal{E}_N)$ is a Frobenius structure on
$N$, i.e., $N$ is a natural Frobenius submanifold of $M$.
\end{proof}

\begin{remark}
\noindent (1) $N$ is the submanifold of $(M,\overline{g}, \circ,e,
\mathcal{E})$. If we assume that $e$ and $\mathcal{E}$ are tangent
to $N$ in Theorem \ref{theorem A}, we recover Theorem
\ref{Theorem4}, then N is a natural Frobenius submanifold.

\noindent (2) In the proof of Theorem \ref{theorem A}, we deduce
that any two equalities can imply the third one:
$$\overline \nabla\mathcal{E}+(\overline \nabla\mathcal{E})^*=D \cdot \Id.$$
$$\nabla\mathcal{E}_N+(\nabla\mathcal{E}_N)^*=D_N \cdot \Id.$$
$$A_{\mathcal{E}_N^\bot}=\lambda \cdot \Id.$$
\end{remark}

\begin{proposition}\label{proposition7}
Let $(M,\overline{g},\overline{\nabla},\circ,e, \mathcal{E})$ be a
Frobenius manifold, $N$ is a submanifold of $M$ such that $g$ is
nondegenerate. If
$$TN\circ TN \subseteq TN$$
$$\overline \nabla=\nabla$$
then $N$ is a natural Frobenius submanifold.
\end{proposition}

\begin{proof}[Proof of proposition \ref{proposition7}]
By the condition $$TN\circ TN \subseteq TN$$
we know $e_N$ is the unit vector field of $(N, TN, \circ)$. And by
$$\overline \nabla=\nabla$$ we get $\nabla\nabla\mathcal{E}_{N}=0$
and $\nabla{e_N}=0$. As in the proof of Theorem 1.1, the structure
connection $\widetilde{\nabla}$ of $N$ is integrable and the unit
$e_N$ is $\nabla$-flat.

$M$ is a Frobenius manifold, so there exist a constant $D$ such
that
$\overline{\nabla}\mathcal{E}+(\overline{\nabla}\mathcal{E})^*=D
\cdot \Id$, for any $X,Y\in TN$. By $$\overline \nabla=\nabla$$ we have:
\begin{eqnarray*}
\<\nabla_X \mathcal{E}_N,Y\> + \<\nabla_Y \mathcal{E}_N,X\>
& = & \<\overline\nabla_X \mathcal{E}_N,Y\> + \<\overline\nabla_Y \mathcal{E}_N,X\>\\
& = & D\cdot \<X,Y\>
\end{eqnarray*}
So $$\nabla\mathcal{E}_N+(\nabla\mathcal{E}_N)^*=D \cdot \Id$$
By the equivalence between Saito structure with metric and Frobenius structure on $TN$, we get $(N,g,\nabla,\circ,e_{N}, \mathcal{E}_{N})$ is Frobenius manifold, i.e., $N$ is the natural Frobenius submanifold of $M$.
\end{proof}

\begin{remark}
We can prove Proposition \ref{proposition7} by applying Theorem
\ref{theorem A} because $\overline \nabla=\nabla$ implies
$A_{\mathcal{E}_N^\bot} =0$. $\overline \nabla=\nabla$ is also not
a necessary condition for a submanifold to be a natural Frobenius
submanifold.
\end{remark}

\begin{example}[\cite{Stra1}]\label{example2}
${B_3\longrightarrow I_2(6)}$

The prepotential for the Frobenius manifold constructed from $B_3$
is
\[
F_{B_3} = \frac{1}{2} {t_1^2t_3}+\frac{1}{2} {t_1 t_2^2} + \frac{1}{6} {t_2^3t_3} +
\frac{1}{6} {t_2^2t_3^3} + \frac{1}{210} {t_3^7}\,.
\]

The two dimensional submanifold is given by
\begin{eqnarray*}
t_1 & = & \tau_1 - \frac{2}{3} k_2^2 \tau_2^3\,,\\
t_2 & = & k_2 \tau_2^2\,,\\
t_3 & = & \tau_2\,.
\end{eqnarray*}
The condition required for the submanifold to be a natural
Frobenius submanifold reduce to $k_2(2k_2-3)(-2k_2-1)=0$. Thus
there are three natural Frobenius submanifolds given by
$k_2=0,-1/2,+3/2$. For $k_2=-1/2$ or $k_2=3/2$, the given
natural Frobenius submanifolds are not totally geodesic
submanifolds.
\end{example}

\section{Frobenius hypersurfaces}\label{section2}

In this section, we mainly talk about the classification of the natural Frobenius hypersurfaces.

\renewcommand{\theequation}
{2.\arabic{equation}}
\setcounter{equation}{0}

For general natural Frobenius submanifold neither $\overline
\nabla=\nabla$ nor $\mathcal{E}=\mathcal{E}_N$ is a necessary
condition. But for hypersurfaces, we get that all the natural
Frobenius submanifolds satisfy either $\overline \nabla=\nabla$ or
$\mathcal{E}=\mathcal{E}_N$.

In this section we suppose $(M,\overline{g},\circ,e, \mathcal{E})$
is a Frobenius manifold, $N$ is a hypersurface of $M$ such that
the restricted metric $g$ is non-degenerate, $\nabla$ is the
Levi-Civita connection of $g$.

\begin{lemma}\label{lemma TG}
Let $(M,\overline g)$ be a Riemannian manifold, let $\overline
\nabla$ be the Levi-Civita connection of $\overline g$, and let
$N$ be a hypersurface of $M$. If there exists a $\overline
\nabla$-flat vector field $X \in TM$ such that $X_N:=pr(X\vert_N)$
is $\nabla$-flat and $X$ is not tangent to $N$, then $N$ is a
totally geodesic submanifold of $M$.
\end{lemma}

\begin{proof} From the flatness of $X$ and $X_N$, we get: $A_{X^\bot}=0$. Because the codimension of $N$ is equal to $1$,
the shape operator $A$ vanishes. So $N$ is a totally geodesic
submanifold.
\end{proof}

\begin{lemma}\label{lemma DD}
If $N$ is a Frobenius submanifold of $M$, and $e$ is tangent to
$N$, then $D=D_N$.
\end{lemma}

\begin{proof} If $N$ is a Frobenius submanifold of $M$, and $e$ is tangent to $N$, so from
$$\overline{\nabla}e=0, \nabla{e_N}=0$$
we get
\begin{equation*}
h(X,e_N)=0,
\end{equation*}
for any $X\in TN$.

From $$\overline{\nabla}\mathcal{E}+(\overline{\nabla}\mathcal{E})^*=D\cdot \Id;$$
$$\nabla\mathcal{E}_N+(\nabla\mathcal{E}_N)^*=D_N \cdot \Id,$$ we get:
\begin{equation*}
\<h(X,Y),{\mathcal{E}_N^\bot}\>=- \frac{D-D_N}{2}\cdot \<X,Y\>.
\end{equation*}
Take $X=e_N$ we get
\begin{equation*}
(D-D_N)\cdot \<e_N,Y\>=0,
\end{equation*}
for any $Y \in TN$.

But $g$ is non-degenerate in $TN$, so we can choose a
local vector field $Y_0$ such that $\<e_N,Y_0\>=1$. So $D=D_N$.
\end{proof}

\begin{proof}[Proof of Proposition \ref{proposition B}(a)]
(2)$\Rightarrow$ (1) By Proposition \ref{proposition7}.

(1)$\Rightarrow$ (2) If $N$ is the natural Frobenius submanifold
of $M$, then $TN\circ TN \subseteq TN$. Because $e$ is not tangent
to $N$ and $N$ is hypersurface of $M$, by Lemma 2.1 $N$ is totally
geodesic submanifold, i.e., $\overline \nabla=\nabla$.
\end{proof}

\begin{proof}[Proof of Proposition \ref{proposition B}(b)]

(2)$\Rightarrow$ (1) If $TN\circ TN \subseteq TN$ and $ \overline
\nabla=\nabla$ then by Proposition \ref{proposition7} $N$ is a
natural Frobenius submanifold of $M$; otherwise if $TN\circ TN
\subseteq TN$ and $\mathcal{E}$ is tangent to $N$, by Theorem
\ref{Theorem4}, we also get $N$ is a natural Frobenius
submanifold.

(1)$\Rightarrow$ (2) Suppose $N$ is a natural Frobenius
submanifold. Because $e$ is tangent to $N$, by Lemma \ref{lemma
DD} we know $D=D_N$, and by
$$\overline{\nabla}\mathcal{E}+(\overline{\nabla}\mathcal{E})^*=D\cdot
\Id;$$
$$\nabla\mathcal{E}_N+(\nabla\mathcal{E}_N)^*=D \cdot \Id; $$ we get:
$$\<h(X,Y),\mathcal{E}_N^\bot \>=0,$$ for any $X,Y \in TN$.

Because codimension of $N$ is $1$, $g$ is non-degenerate, and
$h(X,Y),\mathcal{E}_N^\bot \in TN^ \bot$, so either $h(X,Y)=0$ for
any $X,Y$, or $\mathcal{E}_N^\bot=0$. That is to say either $
\overline \nabla=\nabla$ or $\mathcal{E}$ is tangent to $N$.

$N$ is the natural Frobenius submanifold implies $TN\circ TN \subseteq TN$.
\end{proof}
\begin{remark}
Proposition \ref{proposition B} classify the natural Frobenius
hypersurfaces. But this classification can not be generalized to
natural Frobenius submanifolds of any dimension. We will give an
example of a Frobenius submanifold of codimension two such that
$e_N^\bot\neq 0$, $\overline \nabla \neq \nabla$ and
${\mathcal{E}_N^\bot}\neq 0$.
\end{remark}

\begin{example} \label{example2 continiuation}
In example \ref{example2}, we get two natural Frobenius
submanifolds of $B_3$. Take $k_2=-1/2$. This submanifold, denoted
by $N$, is not a totally geodesic submanifold of $B_3$.

Firstly, just as in example \ref{example1}, we construct a
Frobenius manifold $(B_3 \times \mathcal{A},\overline{g},\circ,e,
\mathcal{E})$ such that $t^1, t^2,t^3, z$ is the flat coordinates
of $B_3 \times \mathcal{A}$. Now embedding $B_3$ to $B_3 \times
\mathcal{A}$:
$$\iota:B_3\longrightarrow B_3 \times \mathcal{A}, P\longmapsto (P, 1).$$
consider the image $\iota(N)$ of $N$ as a submanifold of $B_3
\times \mathcal{A}$. It is given by
\begin{eqnarray*}
t_1 & = & \tau_1 - \frac{1}{6} \tau_2^3\,,\\
t_2 & = & -\frac{1}{2} \tau_2^2\,,\\
t_3 & = & \tau_2\\
z & = & 1.
\end{eqnarray*}
For the Frobenius manifold $B_3 \times \mathcal{A}$, we get a
natural Frobenius submanifold $N\times \lbrace 1\rbrace$ with
$e_N^\bot\neq 0$, $\overline \nabla \neq \nabla$ and
${\mathcal{E}_N^\bot}\neq 0$.
\end{example}

For the first case $e$ is not tangent to $N$, we get some
properties about $\mathcal{E} \circ$ and $\overline \nabla
\mathcal{E}$.

\begin{corollary}\label{corollary R}

Let $(M, \overline g, \circ, e, \mathcal{E})$ be a Frobenius
manifold, and let $N$ be a Frobenius hypersurfaces of $M$ such
that the restricted metric $g$ is non-degenerate. If $e$ is not
tangent to $N$, then $(\mathcal{E} \circ)\vert_{TN}, (\overline
\nabla \mathcal{E})\vert_{TN} \in End(TN)$
\end{corollary}

\begin{proof} From proposition \ref{proposition B} we have two
relations:
$$TN\circ TN \subseteq TN,\quad \overline \nabla=\nabla.$$

If $\mathcal {E}$ is tangent to $N$, then obviously, $(\mathcal{E}
\circ)\vert_{TN}, (\overline \nabla \mathcal{E})\vert_{TN} \in
End(TN)$. Now we suppose $\mathcal {E}$ is not tangent to $N$,
then there exist nonzero function $f$ such that $f
e_N^\bot$=$\mathcal{E}_N^\bot$. But $e_N^\bot \circ TN=0$, so
$\mathcal{E}_N^\bot \circ TN=0$, and then $(\mathcal{E}
\circ)\vert_{TN} \in End(TN)$. Now consider $(\overline \nabla
\mathcal{E})\vert_{TN}$. For any $X \in TN$:

From $\overline \nabla=\nabla$ we know that the second fundamental
form and shape operator vanish, i.e., $h=0$, and $A=0$. Moreover,
from $\overline \nabla=\nabla$, $\overline \nabla e = 0$ and
$\nabla e_N =0$ we get:
\begin{equation*}
\overline \nabla e_N^ \bot = 0.
\end{equation*}
$N$ is a hypersurface of $M$, so there exist a function, denoted by $f$, such that:
\begin{equation*}
\mathcal{E}_N^ \bot = f \cdot e_N^ \bot.
\end{equation*}
we will show that the function $f$ is a constant, i.e., there exists a constant $\mu \in \mathbb{C}$ such that.
\begin{equation*}
\mathcal{E}_N^ \bot = \mu \cdot e_N^ \bot.
\end{equation*}

\emph{Claim}: $\overline \nabla_{\mathcal{E}_N^ \bot} \mathcal{E} \in TN^\bot$.

In fact $M$ is a Frobenius manifold, so we have:
\begin{equation*}
(\overline \nabla_e \mathcal{E})|_{N} = e|_{N}
\end{equation*}
By $\overline \nabla=\nabla$ we can simplify this equality and get:
\begin{equation*}
\nabla_{e_N}\mathcal{E}_N + \nabla^{\bot}_{e_N} \mathcal{E}_N^\bot + \overline \nabla_{e_N^\bot} \mathcal{E} = e_N + e_N^{\bot}
\end{equation*}
But $N$ is a Frobenius submanifold of $M$, so we have $\nabla_{e_N}\mathcal{E}_N = e_N$. So we get:
\begin{equation*}
\overline \nabla_{e_N^\bot} \mathcal{E} = e_N^{\bot}-\nabla^{\bot}_{e_N} \mathcal{E}_N^\bot \in TN^\bot.
\end{equation*}
So
\begin{equation*}
\overline \nabla_{\mathcal{E}_N^ \bot} \mathcal{E} = f \overline \nabla_{e_N^ \bot} \mathcal{E} \in TN^\bot.
\end{equation*}
The flatness of the structure connection $\widetilde {\overline \nabla}$ of $M$ implies that:
\begin{equation*}
\overline \nabla (-\mathcal{E} \circ )- [\Phi, \overline \nabla \mathcal{E}] = - \Phi.
\end{equation*}
which applied to the pair of vectors $(X,\mathcal{E}_N^\bot)$ amounts to
\begin{equation*}
\overline \nabla_X(\mathcal{E}_N^\bot \circ \mathcal{E}) - \overline \nabla_X {\mathcal{E}_N^\bot} \circ \mathcal{E}|_{N} + X \circ \overline \nabla _{\mathcal{E}_N^\bot} \mathcal{E} - \overline \nabla_{X \circ \mathcal{E}_N^\bot}\mathcal{E}=X \circ \mathcal{E}_N^\bot.
\end{equation*}
where $X \in TN$.

From the relation $TN \circ TN^\bot = 0$ we get
\begin{equation*}
\overline \nabla_X(\mathcal{E}_N^\bot \circ \mathcal{E}_N^\bot)-\overline \nabla_X(\mathcal{E}_N^\bot) \circ \mathcal{E}_N^\bot + X \circ \overline \nabla_{\mathcal{E}_N^ \bot} \mathcal{E} = 0.
\end{equation*}
We have proved that $\overline \nabla_{\mathcal{E}_N^ \bot} \mathcal{E} \in TN^\bot$, so the above equality can be simplified to be:
\begin{equation*}
\overline \nabla_X(\mathcal{E}_N^\bot \circ \mathcal{E}_N^\bot)=\overline \nabla(\mathcal{E}_N^\bot) \circ \mathcal{E}_N^\bot
\end{equation*}
By $\mathcal{E}_N^ \bot = f \cdot e_N^ \bot$, $\overline \nabla e_N^\bot=0$ and $e_N^ \bot \circ e_N^ \bot = e_N^ \bot$ we get
\begin{equation*}
X(f^2)=0.
\end{equation*}
So there exists a constant $\mu \in \mathbb{C}$ such that $f=\mu$,
i.e., $\overline \nabla \mathcal{E}_N^\bot=\mu \overline \nabla e_N^\bot=0$.

Then we get:
\[(\overline \nabla \mathcal{E})\vert_{TN}=\overline \nabla {\mathcal{E}_N}+\overline \nabla {\mathcal{E}_N^\bot }=\nabla {\mathcal{E}_N} \in End(TN).
\qedhere
\]
\end{proof}

\begin{remark}
We have another way to see $\overline \nabla \mathcal{E}_N^\bot=0$
in the proof of corollary \ref{corollary R}:

$N$ is a hypersurface, and $\mathcal{E}$ is not tangent to $N$, so
we just need to check
\begin{equation}
\< \overline \nabla_X \mathcal{E}_N^\bot, \mathcal{E}_N^\bot \> =0.
\end{equation}
for any $X \in TN$.

From the relation $\overline{\nabla}\mathcal{E}+(\overline{\nabla}\mathcal{E})^*=D\cdot \Id$, we get
\begin{align*}
\< \overline \nabla_X \mathcal{E}_N^\bot, \mathcal{E}_N^\bot \>
&= \< \overline \nabla_X \mathcal{E}, \mathcal{E}_N^\bot \> \\
&= D \<X, \mathcal{E}_N^\bot \> - \< \overline \nabla_{\mathcal{E}_N^\bot} \mathcal{E}, X\> \\
&= - \< \overline \nabla_{\mathcal{E}_N^\bot} \mathcal{E}, X\>
\end{align*}
where $X \in TN$.

In the proof of corollary \ref{corollary R}, we have the relation
\begin{equation*}
\overline \nabla_{\mathcal{E}_N^ \bot} \mathcal{E} \in TN^\bot.
\end{equation*}
So for any $X \in TN$, we get
\begin{equation*}
\< \overline \nabla_X \mathcal{E}_N^\bot, \mathcal{E}_N^\bot \>=0
\end{equation*}
i.e.,
\begin{equation*}
\overline \nabla_X \mathcal{E}_N^\bot=0.
\end{equation*}
\end{remark}

\bibliographystyle{amsplain}
\bibliography{Frobenius_ZR_070805}

\providecommand{\bysame}{\leavevmode\hbox to3em{\hrulefill}\thinspace}
\providecommand{\MR}{\relax\ifhmode\unskip\space\fi MR }
\providecommand{\MRhref}[2]{%
  \href{http://www.ams.org/mathscinet-getitem?mr=#1}{#2}
}
\providecommand{\href}[2]{#2}
\begin{thebibliography}{1}

\bibitem{D}
Boris Dubrovin, \emph{Geometry of {$2$}{D} topological field theories},
  Integrable systems and quantum groups (Montecatini Terme, 1993), Lecture
  Notes in Math., vol. 1620, Springer, Berlin, 1996, pp.~120--348.
  \MR{MR1397274 (97d:58038)}

\bibitem{Hert}
Claus Hertling, \emph{Frobenius manifolds and moduli spaces for singularities},
  Cambridge Tracts in Mathematics, vol. 151, Cambridge University Press,
  Cambridge, 2002. \MR{MR1924259 (2004a:32043)}

\bibitem{Mani}
Yuri~I. Manin, \emph{Frobenius manifolds, quantum cohomology, and moduli
  spaces}, American Mathematical Society Colloquium Publications, vol.~47,
  American Mathematical Society, Providence, RI, 1999. \MR{MR1702284
  (2001g:53156)}

\bibitem{Sabb}
Claude Sabbah, \emph{D\'eformations isomonodromiques et vari\'et\'es de
  {F}robenius}, Savoirs Actuels (Les Ulis). [Current Scholarship (Les Ulis)],
  EDP Sciences, Les Ulis, 2002, , Math\'ematiques (Les Ulis). [Mathematics (Les
  Ulis)]. \MR{MR1933784 (2003m:32013)}

\bibitem{Stra1}
I.~A.~B. Strachan, \emph{Frobenius submanifolds}, J. Geom. Phys. \textbf{38}
  (2001), no.~3-4, 285--307. \MR{MR1829045 (2002f:53150)}

\bibitem{Stra2}
\bysame, \emph{Frobenius manifolds: natural submanifolds and induced
  bi-{H}amiltonian structures}, Differential Geom. Appl. \textbf{20} (2004),
  no.~1, 67--99. \MR{MR2030167 (2004m:53156)}

\bibitem{Will}
T.~J. Willmore, \emph{Riemannian geometry}, Oxford Science Publications, The
  Clarendon Press Oxford University Press, New York, 1993. \MR{MR1261641
  (95e:53002)}

\end{thebibliography}
\end{document}